# NUMERICAL APPROXIMATIONS FOR THE CAHN-HILLIARD PHASE FIELD MODEL OF THE BINARY FLUID-SURFACTANT SYSTEM

XIAOFENG YANG

ABSTRACT. In this paper, we consider the numerical approximations for the commonly used binary fluid-surfactant phase field model that consists two nonlinearly coupled Cahn-Hilliard equations. The main challenge in solving the system numerically is how to develop easy-to-implement time stepping schemes while preserving the unconditional energy stability. We solve this issue by developing two linear and decoupled, first order and a second order time-stepping schemes using the so-called "Invariant Energy Quadratization" approach for the double well potentials and a subtle explicit-implicit technique for the nonlinear coupling potential. Moreover, the resulting linear system is well-posed and the linear operator is symmetric positive definite. We rigorously prove the first order scheme is unconditionally energy stable. Various numerical simulations are presented to demonstrate the stability and the accuracy thereafter.

## 1. INTRODUCTION

Surfactants are usually organic compounds that can alter or reduce the surface tension of the solution, and allows for the mixing of dissimilar (immiscible) liquids. A typical well-known example of dissimilar liquids is the mixture of oil and water, where the water molecule is polar, and it does not hang out with nonpolar molecules like oil. Hence, in order to make the mixing more favorable, a molecular intermediate, commonly referred as surfactants, is needed. The surfactant molecules are amphiphilic, they are capable of binding to both water and oil molecules (hydrophilic heads into the water, hydrophobic tails into the oil), allowing these molecules to mix via reducing the surface tension between water and oil. This is the essential reason that people can clean the hands or clothes using soaps or detergents. Except the daily uses in cleaning, surfactants have been widely used in industrial fields for various purposes, e.g., oil recovery [52] and food processing [36], etc.

There is a large quantity of studies on the modeling and numerical simulations to investigate the binary fluid-surfactant system. In the pioneering work of Laradji et. al. [25, 26], the diffuse interface approach, or called the phase field method, was firstly used to study the phase transition behaviors of the monolayer microemulsion system, formed by surfactant molecules. Following Laradji's work, a number of phase field type of fluid-surfactant models have been developed during the last two decades, see also [11, 13, 24, 53–55]. In this paper, we consider the numerical approximations for solving the model that was developed by Komura et. al. in [24] since it appears to be the minimal model with the least number of nonlinear coupling entropy terms in the free energy. The main difference between the surfactant models for fluid-surfactant system and the most classical Cahn-Hilliard phase field model for two phase system, that was proposed by Cahn and Hilliard in [4] and well-studied in [5,6,8,12,20,22,27,28,32,34,51,61,64,65] and the references therein, is that an extra phase variable (surfactant) has to be used. Hence two physical phase field variables are incorporated in the model, where one is used to represent the local densities of the fluids, and the other is used for the local concentration of the surfactant. Except the regular potentials, e.g., the hydrophilic (gradient entropy) and hydrophobic part (nonlinear double well potential) for each phase field variable, some nonlinear coupling entropy terms are participated in the total free energy. By minimizing the total energy via the variational

Department of Mathematics, University of South Carolina, Columbia, SC 29208, USA. Email: xfyang@math.sc.edu. X. Yang's research is partially supported by the U.S. National Science Foundation under grant numbers DMS-1200487 and DMS-1418898.

approach, one can obtain the governing system that consists *two nonlinearly coupled Cahn-Hilliard* type equations.

From the numerical point of view, although the phase field variable is continuous and smooth, the induced models are still very stiff where the stiffness is induced by an order parameter that represents the thickness of the interface. It can be seen clearly from the a fact that some severe stability conditions on the time step will occur if the nonlinear term is discretized in some normal ways like fully implicit or explicit type approaches. Such a constraint on the time step can cause very high computational cost in practice [10, 44]. Hence, it is desirable to establish efficient numerical schemes that can verify the so-called energy stable property at the discrete level irrespectively of the coarseness of the discretization. Such kinds of algorithms are usually called unconditionally energy stable or thermodynamically consistent. The scheme with this property is specially preferred since it is not only critical for the numerical scheme to capture the correct long time dynamics of the system, but also provides sufficient flexibility for dealing with the stiffness issue. Although a variety of the phase field models of the fluid-surfactant system [24–26, 54] had been built for over twenty years, we must notice that there are very few successful attempts in designing unconditional energy stable schemes where the main difficulty lies in the discretization for the nonlinear terms, in particular the nonlinear coupling term associated with multiple phase field variables. In [21], a second order in time, nonlinear scheme was developed based on the Crank-Nicolson method. However, the scheme could be computationally expensive due to its highly nonlinear and coupling nature. Moreover, since the nonlinear terms are mostly treated implicitly, the energy stability of the scheme does not hold at all, and the solvability is unknown yet.

Therefore, the main purpose of this paper is to construct some time discretization schemes that are expected to be *easy-to-implement* (linear system) and *unconditionally energy stable* (with a discrete energy dissipation law). We achieve this goal by adopting the "Invariant Energy Quadratization" (IEQ) approach (cf. [17, 58, 62, 63, 66, 71, 72]) for the nonlinear double well potentials, and a subtle explicit-implicit treatment for the nonlinear coupling term. As a result, in each time step, the scheme leads to two *decoupled* linear systems. Moreover, we show that each corresponding linear operator is *symmetric positive definite*, so that one can solve it using the well-developed fast matrix solvers efficiently (CG or other Krylov subspace methods). To the best of the author's knowledge, the numerical scheme proposed in this paper is the first linear, decoupled schemes for the nonlinearly coupled multivariable phase-field model.

The rest of the paper is organized as follows. In Section 2, we present the whole model and give the PDE energy law. In Section 3, we develop a linear, totally decoupled, first order time stepping scheme and prove the well-posedness and unconditional stability. We further develop a linear, totally decoupled second order time stepping schemes based on the Adam-Bashforth method. Various numerical experiments are carried out in Section 4 to validate the accuracy and stability of the proposed scheme. Finally, some concluding remarks are presented in Section 5.

## 2. Models

We now give a brief introduction for binary fluid-surfactant phase field model. To fix the notations, without ambiguity, we denote by $(f(\boldsymbol{x}), g(\boldsymbol{x})) = (\int_\Omega f(\boldsymbol{x})g(\boldsymbol{x})d\boldsymbol{x})^{\frac{1}{2}}$ the $L^2$ inner product between functions $f(\boldsymbol{x})$ and $g(\boldsymbol{x})$, by $\|f\| = (f,f)$ the $L^2$ norm of function $f(\boldsymbol{x})$, where $\Omega$ is the computed domain.

In the ternary fluid flow system for the mixture of water, oil and surfactant, monolayers of surfactant molecules form microemulsions as a random phase. Such a microemulsion system usually exhibits various interesting microstructures depending on the temperature or the composition. In [24], the dynamics of microphase separation in microemulsion systems was modeled by a phase field model with two order parameters, where, the free energy of the system is given as follows,

$$(2.1) \quad E(\phi, \rho) = \int_\Omega \Big( \underbrace{\frac{1}{2}|\nabla\phi|^2 + \frac{\alpha}{2}(\Delta\phi)^2 + \frac{1}{4\epsilon^2}(\phi^2-1)^2}_{\text{part A}} + \underbrace{\frac{\beta}{2}|\nabla\rho|^2 + \frac{1}{4\eta^2}\rho^2(\rho-\rho_s)^2}_{\text{part B}} - \underbrace{\theta\rho|\nabla\phi|^2}_{\text{part C}} \Big) d\boldsymbol{x}$$

where $\alpha, \beta, \epsilon, \eta, \rho_s, \theta$ are all positive parameters.

Now we give some brief descriptions about the free energy. More detailed descriptions about the modeling can be found in [24].

In part A, the phase field variable $\phi$ is introduced to label the local densities of the two fluids, e.g., water and oil, where the fourth order polynomial is the double well Ginzburg-Landau potential. Thus there are two minimum values for this polynomial which means the two bulk state for $\phi$, i.e.,

$$\phi(\boldsymbol{x}, t) = \begin{cases} 1 & \text{fluid I,} \\ -1 & \text{fluid II.} \end{cases} \tag{2.2}$$

The interface of the fluid mixture is described by the zero level set $\Gamma_t = \{\boldsymbol{x} : \phi(\boldsymbol{x}, t) = 0\}$. Part I is the well-known "mixing" energy for phase field model, where the linear part $((\Delta\phi)^2$ and $|\nabla\phi|^2)$ contributes to the hydrophilic type (tendency of mixing) of interactions between the materials and the double well potential represents the hydrophobic type (tendency of separation) of interactions. As the consequence of the competition between the two types of interactions, the equilibrium configuration will include a diffusive interface. About the theoretical or numerical study about the Cahn-Hilliard or Allen-Cahn system derived from this part of energy, we refer to [1, 6, 16, 19, 28–30, 37, 40, 42–47, 50, 59, 60, 69, 70, 73, 74].

In part B, the second phase field variable $\rho(\boldsymbol{x}, t)$ is introduced to represent the local concentration of surfactants where another fourth double well polynomial potential is used as well. There are two minimum values that enable two bulk states for $\rho$, i.e., $\rho = 0$ and $\rho = \rho_s$. The state of $\rho = 0$ corresponds to the case in which the system is locally occupied either by oil of water without surfactants. The state of $\rho = \rho_s$ corresponds to the case in which the local volume is occupied only by surfactants. The quantity $\rho_s$ can be considered to represent the density of condensed hydrocarbon chains of surfactants when they self-assemble. In the numerical experiments, we simply set $\rho_s = 1$.

Part C represents the coupling between the surfactants and fluid interface, i.e., the state $\rho = \rho_s$ tends to occupy the narrow region around the oil-water interfaces. It is also essential in microemulsions that the interfacial tension vanishes when the interface is saturated with surfactants.

The time evolution equation is assumed to be the Cahn-Hilliard type, i.e., the gradient flow in $H^{-1}$, the system reads as follows,

$$\phi_t = M_\phi \Delta \mu_\phi, \tag{2.3}$$

$$\mu_\phi = \frac{\delta E}{\delta \phi} = -\Delta\phi + \alpha\Delta^2\phi + \frac{1}{\epsilon^2}\phi(\phi^2 - 1) + 2\theta\nabla\cdot(\rho\nabla\phi), \tag{2.4}$$

$$\rho_t = M_\rho \Delta \mu_\rho, \tag{2.5}$$

$$\mu_\rho = \frac{\delta E}{\delta \rho} = -\beta\Delta\rho + \frac{1}{\eta^2}\rho(\rho - \rho_s)(\rho - \frac{\rho_s}{2}) - \theta|\nabla\phi|^2. \tag{2.6}$$

For the domain $\Omega$, the boundary conditions can be

$$(i) \text{ all variables are periodic; or } (ii) \partial_{\boldsymbol{n}}\phi|_{\partial\Omega} = \partial_{\boldsymbol{n}}\rho|_{\partial\Omega} = \nabla\mu_\phi \cdot \boldsymbol{n}|_{\partial\Omega} = \nabla\mu_\rho \cdot \boldsymbol{n}|_{\partial\Omega} = 0, \tag{2.7}$$

where $\boldsymbol{n}$ is the outward normal on the domain boundary $\partial\Omega$. The system (2.3)-(2.5) conserves the local mass density, i.e., $\frac{d}{dt}\int_\Omega \phi d\boldsymbol{x} = \frac{d}{dt}\int_\Omega \rho d\boldsymbol{x} = 0$.

It is straight forward to obtain the PDE energy law for the Cahn-Hilliard system (2.3)-(2.6). By taking the $L^2$ inner product of (2.3) with $\mu_\phi$, of (2.4) with $\phi_t$ of (2.5) with $\mu_\rho$, of (2.6) with $\rho_t$, performing integration by parts, and combining all equalities, we can derive

$$\frac{d}{dt}E(\phi, \rho) = -M_\phi\|\nabla\mu_\phi\|^2 - M_\rho\|\nabla\mu_\rho\|^2 \leq 0. \tag{2.8}$$

**Remark 2.1.** *Alternatively, one can also assume the Allen-Cahn dynamics for the phase variables, namely*

$$\phi_t = -M_\phi\mu_\phi, \tag{2.9}$$
$$\rho_t = -M_\rho\mu_\rho,$$

*with periodic boundary conditions or no-flux boundary condition as $\partial_{\boldsymbol{n}}\phi|_{\partial\Omega} = \partial_{\boldsymbol{n}}\rho|_{\partial\Omega} = 0$. Note that the volume conservation property does not hold for Allen-Cahn type dynamics. To overcome it, one can equip (2.9) with a scalar Lagrange multiplier to enforce this conservation property (cf. [42, 59]), or modify the free energy functional (2.1) by adding a penalty term for mass conservation, similar as the method used in [8, 33]. Both ways will not introduce any mathematical or numerical difficulties and all analysis for Cahn-Hilliard*

*system can be carried out for Allen-Cahn system without further difficulties, thus we shall not study the Allen-Cahn system in this paper.*

**Remark 2.2.** *The total free energy (2.1) is the rescale of the free energy proposed in [24]. We notice that the gradient term $|\nabla \rho|^2$ is actually neglected in [24]. The absence of this term is "physically" reasonable since the energy cost due to the direct attachment between hydrocarbon chains and oil molecules or between hydrophilic head and water molecules is small. But in this paper, we add this term back "mathematically" because this term may help the total energy $E(\phi, \rho)$ to be bounded from below. However, it is not a trivial work to prove the total free energy to be bounded from below even with this gradient potential, where the main difficulty is how to bound the nonlinear coupling potential (part C). A possible way is to follow Caffarelli and Muler's work for the classical Cahn-Hilliard equation [3], namely, to prove the concentration variable $\rho$ is $L^\infty$ bounded. Then the part C can be bounded by other terms. The related PDE analysis is out of the scope of this paper, hence we only present the evolution of the free energy through the numerical experiments in Section 4. Indeed, the coupling term between $\phi$ and $\rho$ in part C can cause serious hurdles to obtain the error estimates for the schemes even though its energy stability can be formally derived. We will consider this issue in the future work.*

**Remark 2.3.** *When concerning the hydrodynamics effects, by assuming both fluids are incompressible, the system reads as follows, see also [38],*

$$\phi_t + \nabla \cdot (\boldsymbol{u}\phi) = M_\phi \Delta \mu_\phi, \tag{2.10}$$

$$\rho_t + \nabla \cdot (\boldsymbol{u}\rho) = M_\rho \Delta \mu_\rho, \tag{2.11}$$

$$\boldsymbol{u}_t + (\boldsymbol{u} \cdot \nabla)\boldsymbol{u} + \nabla p - \nu \Delta \boldsymbol{u} + \phi \nabla \mu_\phi + \rho \nabla \mu_\rho = 0, \tag{2.12}$$

$$\nabla \cdot \boldsymbol{u} = 0, \tag{2.13}$$

*where $\boldsymbol{u}$ is the fluid velocity, $p$ is the pressure, $\nu$ is the viscosity, the two nonlinear terms related to $\phi$ and $\rho$ in the momentum equation are the induced stresses due to the total free energy. The numerical method for the hydrodynamics coupled phase field models in various situations, e.g., multiphase complex fluids, liquid crystals, had been well-studied in literatures, see also [23, 28, 31, 35, 39, 43, 45, 46, 68, 73].*

## 3. Numerical Schemes

We now construct time stepping schemes to solve the model system (2.3)-(2.6). *The aim is to construct schemes that are not only easy-to-implement, but also unconditionally energy stable.* Here the term "easy-to-implement" is referred to "linear" and "decoupled" in comparison with its counter parts: "nonlinear" and "coupled".

We notice that there are two numerical challenges, including how to decouple the computations of $\phi$ and $\rho$; and how to discretize the two double well polynomial terms for $\phi$ and $\rho$. For the nonlinear coupled potential of $\phi$ and $\rho$, we have to develop a subtle explicit-implicit technique to obtain a decoupling and linear scheme. For the double well polynomial potential terms, we recall that there are two commonly used techniques to discretize this term while preserving the energy stability. The first well-known technique is the so-called *convex splitting* approach [9, 17, 18, 41, 73, 74], where the convex part of the potential is treated implicitly and the concave part is treated explicitly. The convex splitting approach is energy stable, however, it produces a nonlinear scheme at most cases. Thus its implementation is complicated and the computational cost is high, that excludes this method in this paper. The second technique is the *linear stabilization* approach (cf. [6, 30, 33, 42–44, 46, 48–50, 56, 57, 59, 60, 67, 69, 70]), where the nonlinear term is simply treated explicitly and some linear stabilizing terms are added to maintain the stability. This approach is linear so it is very efficient and easy to implement. But its stability requests a special property (generalized maximum principle) satisfied by the classical solution of the PDE and the numerical solution, that is usually very hard to prove.

We now describe our special treatments to discretize the nonlinear terms to preserve the energy stability. For the polynomial terms, we adopt the IEQ approach, which has been successfully applied to solve a variety of gradient flow type models (cf. [58, 63, 71, 72]). Its idea, that is very simple but quite different from those traditional methods like implicit, explicit, nonlinear splitting, or other various tricky Taylor expansions to discretize the nonlinear potentials, is to make the free energy *quadratic*. More precisely, we introduce two auxiliary functions as

$$U = \phi^2 - 1, \tag{3.1}$$

$$V = \rho(\rho - \rho_s). \tag{3.2}$$

In turn, the total energy can be rewritten as

$$E(\phi, \rho, U, V) = \int_\Omega \Big(\frac{1}{2}|\nabla\phi|^2 + \frac{\alpha}{2}(\Delta\phi)^2 + \frac{1}{4\epsilon^2}U^2 + \frac{\beta}{2}|\nabla\rho|^2 + \frac{1}{4\eta^2}V^2 - \theta\rho|\nabla\phi|^2\Big)d\boldsymbol{x}. \tag{3.3}$$

Now, we have a new, but equivalent PDE system as follows,

$$\phi_t = M_\phi \Delta \mu_\phi, \tag{3.4}$$

$$\mu_\phi = -\Delta\phi + \alpha\Delta^2\phi + \frac{1}{\epsilon^2}HU + 2\theta\nabla\cdot(\rho\nabla\phi), \tag{3.5}$$

$$\rho_t = M_\rho \Delta \mu_\rho, \tag{3.6}$$

$$\mu_\rho = -\beta\Delta\rho + \frac{1}{\eta^2}GV - \theta|\nabla\phi|^2, \tag{3.7}$$

$$U_t = 2H\phi_t, \tag{3.8}$$

$$V_t = 2G\rho_t, \tag{3.9}$$

where $H(\phi) = \phi, G(\rho) = \rho - \frac{\rho_s}{2}$.

The boundary conditions for the new system are still (2.7) since the equations (3.8) and (3.9) for the new variable $U$ and $V$ are simply ODEs with time. The initial conditions read as

$$\phi|_{(t=0)} = \phi_0, \ \rho|_{(t=0)} = \rho_0, \ U|_{(t=0)} = \phi_0^2 - 1, \ V|_{(t=0)} = \rho_0(\rho_0 - \rho_s). \tag{3.10}$$

To derive the energy dissipative law for this system, we take the $L^2$ inner product of (3.4) with $\mu_\phi$, of (3.5) with $-\phi_t$, of (3.6) with $\mu_\rho$, of (3.7) with $-\rho_t$, of (3.8) with $\frac{1}{2\epsilon^2}U$, of (3.9) with $\frac{1}{2\eta^2}V$ and combining all terms together, we obtain the new energy dissipation law as follows,

$$\frac{d}{dt}E(\phi, \rho, U, V) = -M_\phi\|\nabla\mu_\phi\|^2 - M_\rho\|\nabla\mu_\rho\|^2 \leq 0. \tag{3.11}$$

**Remark 3.1.** *We consider the two quartic, double well polynomial potential for $\phi$ and $\rho$ as two quadratic functionals by applying appropriate substitutions if needed. Therefore, after simple substitutions using new variables $U, V$, the energy is transformed to an equivalent quadratic form. We emphasize that the new transformed system (3.4)-(3.9) is exactly equivalent to the original system (2.3)-(2.6) since (3.1)-(3.2) can be easily obtained by integrating (3.8)-(3.9) with respect to the time. Therefore, the energy law (3.11) for the transformed system is exactly the same as the energy law (2.8) for the original system for the time-continuous case. We will develop time-marching schemes for the new transformed system (3.4)-(3.9) that in turn follows the new energy dissipation law (3.11) instead of the energy law (2.8) for the original system.*

3.1. **First order scheme.** We now present a first order time marching scheme to solve the system (3.4)-(3.9).

<u>Scheme</u> **1.** *Assuming that $\phi^n, \rho^n, U^n, V^n$ are already known, we compute $\phi^{n+1}, \rho^{n+1}, U^{n+1}, V^{n+1}$ from the following first order temporal semi-discrete system:*
    **Step 1:** *We update $\rho^{n+1}$ and $V^{n+1}$ as follows,*

$$\frac{\rho^{n+1} - \rho^n}{\delta t} = M_\rho \Delta \mu_\rho^{n+1}, \tag{3.12}$$

$$\mu_\rho^{n+1} = -\beta\Delta\rho^{n+1} + \frac{1}{\eta^2}G^n V^{n+1} - \theta|\nabla\phi^n|^2, \tag{3.13}$$

$$V^{n+1} - V^n = 2G^n(\rho^{n+1} - \rho^n). \tag{3.14}$$

**Step 2:** We update $\phi^{n+1}$ and $U^{n+1}$ as follows,

$$\frac{\phi^{n+1} - \phi^n}{\delta t} = M_\phi \Delta \mu_\phi^{n+1}, \tag{3.15}$$

$$\mu_\phi^{n+1} = -\Delta \phi^{n+1} + \alpha \Delta^2 \phi^{n+1} + \frac{1}{\epsilon^2} H^n U^{n+1} + 2\theta \nabla \cdot (\rho^{n+1} \nabla \frac{\phi^{n+1} + \phi^n}{2}), \tag{3.16}$$

$$U^{n+1} - U^n = 2H^n(\phi^{n+1} - \phi^n). \tag{3.17}$$

The boundary conditions can be

(i) all variables are periodic; or (ii) $\partial_{\boldsymbol{n}} \phi^{n+1}|_{\partial \Omega} = \partial_{\boldsymbol{n}} \rho^{n+1}|_{\partial \Omega} = \nabla \mu_\phi^{n+1} \cdot \boldsymbol{n}|_{\partial \Omega} = \nabla \mu_\rho^{n+1} \cdot \boldsymbol{n}|_{\partial \Omega} = 0$.

**Remark 3.2.** When computing $\rho^{n+1}$ in step 1, we only need $\phi^n$. When computing $\phi^{n+1}$ in step 2, $\rho^{n+1}$ is already obtained from step 1. Thus the computations of $\phi$ and $\rho$ are totally decoupled. Furthermore, when computing $\rho^{n+1}$ and $\phi^{n+1}$, one does not need calculate $U^{n+1}$ and $V^{n+1}$. Note the coefficient $H$ and $G$ of the new variables $U$ and $V$ are both treated explicitly, thus we can rewrite the equations (3.14) and (3.17) as follows:

$$\begin{cases} V^{n+1} = 2G^n \rho^{n+1} + D^n, \\ U^{n+1} = 2H^n \phi^{n+1} + S^n, \end{cases} \tag{3.18}$$

where $D^n = V^n - 2G^n \rho^n$ and $S^n = U^n - 2H^n \phi^n$. Thus (3.12)-(3.13) can be rewritten as the following linear system

$$\frac{1}{\delta t} \rho^{n+1} = M_\rho \Delta \mu_\rho^{n+1} + \frac{1}{\delta t} \rho^n, \tag{3.19}$$

$$\mu_\rho^{n+1} = P(\rho^{n+1}) + g_1^n, \tag{3.20}$$

where $P(\rho) = -\beta \Delta \rho + \frac{2}{\eta^2}(G^n)^2 \rho$, $g_1^n = \frac{1}{\eta^2} G^n D^n - \theta |\nabla \phi^n|^2$. Similarly, (3.15)-(3.16) can be rewritten as the following linear system

$$\frac{1}{\delta t} \phi^{n+1} = M_\phi \Delta \mu_\phi^{n+1} + \frac{1}{\delta t} \phi^n, \tag{3.21}$$

$$\mu_\phi^{n+1} = Q(\phi^{n+1}) + g_2^n, \tag{3.22}$$

where $Q(\phi) = -\Delta \phi + \alpha \Delta^2 \phi + \frac{2}{\epsilon^2}(H^n)^2 \phi + \theta \nabla \cdot (\rho^{n+1} \nabla \phi)$ and $g_2^n = \frac{1}{\epsilon^2} H^n S^n + \theta \nabla \cdot (\rho^{n+1} \nabla \phi^n)$. Therefore, we can solve $\rho^{n+1}$ and $\psi^{n+1}$ directly from (3.19)-(3.20) and (3.21)-(3.22). Once we obtain $\phi^{n+1}, \rho^{n+1}$, then $V^{n+1}, U^{n+1}$ are automatically given in (3.18). Namely, the new variables $V$ an $U$ will not involve any extra computational costs.

We first show the well-posedness of the linear system (3.12)-(3.14) (or (3.19)-(3.20)) as follows.

**Theorem 3.1.** The linear system (3.19)-(3.20) admits a unique solution in $H^1(\Omega)$, and the linear operator is symmetric positive definite.

*Proof.* From (3.19), by taking the $L^2$ inner product with 1, we have

$$\int_\Omega \rho^{n+1} d\boldsymbol{x} = \int_\Omega \rho^n d\boldsymbol{x} = \cdots = \int_\Omega \rho^0 d\boldsymbol{x}. \tag{3.23}$$

Let $V_\rho = \frac{1}{|\Omega|} \int_\Omega \rho^0 d\boldsymbol{x}$, $V_\mu = \frac{1}{|\Omega|} \int_\Omega \mu_\rho^{n+1} d\boldsymbol{x}$, and we define

$$\widehat{\rho}^{n+1} = \rho^{n+1} - V_\rho, \ \widehat{\mu}_\rho^{n+1} = \mu_\rho^{n+1} - V_\mu. \tag{3.24}$$

Thus, from (3.19)-(3.20), $(\widehat{\rho}^{n+1}, \widehat{\mu}_\rho^{n+1})$ are the solutions for the following equations with unknowns $(\rho, w)$,

$$\frac{1}{M_\rho \delta t} \rho - \Delta \mu = f, \tag{3.25}$$

$$\mu + V_\mu - P(\rho) = g, \tag{3.26}$$

where $f = \frac{1}{M_\rho \delta t}\hat{\rho}^n$, $g = g_1^n + \frac{2}{\eta^2}G^nG^nV_\rho$. Moreover, $f, \rho$ and $\mu$ are all mean 0, i.e., $\int_\Omega f d\boldsymbol{x} = \int_\Omega \rho d\boldsymbol{x} = \int_\Omega \mu d\boldsymbol{x} = 0$.

Define the inverse Laplace operator $u$ (with $\int_\Omega u d\boldsymbol{x} = 0$) $\mapsto v := \Delta^{-1}u$ by

$$(3.27) \quad \begin{cases} \Delta v = u, \quad \int_\Omega v d\boldsymbol{x} = 0, \\ \text{with the boundary conditions either (i) } v \text{ is periodic, or (ii) } \partial_{\boldsymbol{n}} v|_{\partial\Omega} = 0. \end{cases}$$

Applying $-\Delta^{-1}$ to (3.25) and using (3.26), we obtain

$$(3.28) \quad -\frac{1}{M_\rho \delta t}\Delta^{-1}\rho + P(\rho) - V_\mu = -\Delta^{-1}f - g,$$

then we can express the above linear system (3.28) as $\mathbb{A}\rho = b$.

(i). For any $\rho_1$ and $\rho_2$ in $H^1(\Omega)$ satisfy the boundary conditions (2.7) and $\int_\Omega \rho_1 dx = \int_\Omega \rho_2 dx = 0$, using integration by parts, we derive

$$(3.29) \quad \begin{aligned} (\mathbb{A}(\rho_1), \rho_2) &= -\frac{1}{M_\rho \delta t}(\Delta^{-1}\rho_1, \rho_2) + (P(\rho_1), \rho_2) \\ &\leq C_1(\|\nabla\Delta^{-1}\rho_1\|\|\nabla\Delta^{-1}\rho_2\| + \|\nabla\rho_1\|\|\nabla\rho_2\| + \|\rho_1\|\|\rho_2\|) \\ &\leq C_2\|\rho_1\|_{H^1}\|\rho_2\|_{H^1}. \end{aligned}$$

Therefore, the bilinear form $(\mathbb{A}(\rho_1), \rho_2)$ is bounded for any $\rho_1, \rho_2 \in H^1(\Omega)$ with mean 0.

(ii). For any $\phi \in H^1(\Omega)$ with mean 0, it is easy to derive that, ,

$$(3.30) \quad (\mathbb{A}(\rho), \rho) = \frac{1}{M_\rho \delta t}\|\nabla\Delta^{-1}\rho\|^2 + \beta\|\nabla\rho\|^2 + \frac{2}{\eta^2}\|G^n\rho\|^2 \geq C_3\|\rho\|_{H^1}^2,$$

from Poincare inequality. Thus the bilinear form $(\mathbb{A}(\rho_1), \rho_2)$ is coercive.

Then from the Lax-Milgram theorem, we conclude the linear system (3.28) admits a unique solution in $H^1(\Omega)$.

For any $\rho_1, \rho_2$ with $\int_\Omega \rho_1 d\boldsymbol{x} = 0$ and $\int_\Omega \rho_2 d\boldsymbol{x} = 0$, we can easily derive

$$(3.31) \quad (\mathbb{A}\rho_1, \rho_2) = (\rho_1, \mathbb{A}\rho_2).$$

Thus $\mathbb{A}$ is self-adjoint. Meanwhile, from (3.30), we derive $(\mathbb{A}\rho, \rho) \geq 0$ for any $\rho \in H^1(\Omega)$ and mean 0, where "=" is valid if only if $\rho = 0$. This concludes the linear operator $\mathbb{A}$ is positive definite. □

Likewise, the linear system (3.15)-(3.17) (or (3.21)-(3.22)) is well-posed, symmetric positive definite, and admits a unique solution in $H^2(\Omega)$. The proof is very similar to Theorem 3.1, thus we omit the details here.

We show the unconditional energy stability of the scheme (3.12)-(3.17) as follows.

**Theorem 3.2.** *The scheme (3.12)-(3.17) is unconditionally energy stable and satisfies the following discrete energy dissipation law,*

$$(3.32) \quad \frac{1}{\delta t}(\mathbb{E}^{n+1} - \mathbb{E}^n) \leq -M_\phi\|\nabla\mu_\phi^{n+1}\|^2 - M_\rho\|\nabla\mu_\rho^{n+1}\|^2,$$

*where*

$$(3.33) \quad \mathbb{E}^n = \int_\Omega \Big(\frac{1}{2}|\nabla\phi^n|^2 + \frac{\alpha}{2}(\Delta\phi^n)^2 + \frac{1}{4\epsilon^2}(U^n)^2 + \frac{\beta}{2}|\nabla\rho^n|^2 + \frac{1}{4\eta^2}(V^n)^2 - \theta\rho^n|\nabla\phi^n|^2\Big)d\boldsymbol{x}.$$

*Proof.* By taking the $L^2$ inner product of (3.12) with $\delta t \mu_\rho^{n+1}$, we obtain

$$(3.34) \quad (\rho^{n+1} - \rho^n, \mu_\rho^{n+1}) = -M_\rho \delta t\|\nabla\mu_\rho^{n+1}\|^2.$$

By taking the $L^2$ inner product of (3.13) with $-(\rho^{n+1} - \rho^n)$, and applying the following identities

$$(3.35) \quad 2(a - b, a) = |a|^2 - |b|^2 + |a - b|^2,$$

we obtain

$$
\begin{aligned}
(3.36) \quad -(\mu_\rho^{n+1}, \rho^{n+1} - \rho^n) = &-\frac{\beta}{2}(\|\nabla\rho^{n+1}\|^2 - \|\nabla\rho^n\|^2 + \|\nabla(\rho^{n+1} - \rho^n)\|^2) \\
&-\frac{1}{\eta^2}(G^n V^{n+1}, \rho^{n+1} - \rho^n) + \theta(|\nabla\phi^n|^2, \rho^{n+1} - \rho^n).
\end{aligned}
$$

By taking the $L^2$ inner product of (3.14) with $\frac{1}{2\eta^2}V^{n+1}$, we obtain

$$
(3.37) \quad \frac{1}{4\eta^2}(\|V^{n+1}\|^2 - \|V^n\|^2 + \|V^{n+1} - V^n\|^2) = \frac{1}{\eta^2}(G^n(\rho^{n+1} - \rho^n), V^{n+1}).
$$

By taking the $L^2$ inner product of (3.15) with $\delta t \mu_\phi^{n+1}$, we obtain

$$
(3.38) \quad (\phi^{n+1} - \phi^n, \mu_\phi^{n+1}) = -M_\phi \delta t \|\nabla\mu_\phi^{n+1}\|^2.
$$

By taking the $L^2$ inner product of (3.16) with $-(\phi^{n+1} - \phi^n)$, we obtain

$$
\begin{aligned}
(3.39) \quad -(\mu_\phi^{n+1}, \phi^{n+1} - \phi^n) = &-\frac{1}{2}(\|\nabla\phi^{n+1}\|^2 - \|\nabla\phi^n\|^2 + \|\nabla(\phi^{n+1} - \phi^n)\|^2) \\
&-\frac{\alpha}{2}(\|\Delta\phi^{n+1}\|^2 - \|\Delta\phi^n\|^2 + \|\Delta(\phi^{n+1} - \phi^n)\|^2) \\
&-\frac{1}{\epsilon^2}(H^n U^{n+1}, \phi^{n+1} - \phi^n) + \theta(\rho^{n+1}, |\nabla\phi^{n+1}|^2 - |\nabla\phi^n|^2).
\end{aligned}
$$

By taking the $L^2$ inner product of (3.17) with $\frac{1}{2\epsilon^2}U^{n+1}$, we obtain

$$
(3.40) \quad \frac{1}{4\epsilon^2}(\|U^{n+1}\|^2 - \|U^n\|^2 + \|U^{n+1} - U^n\|^2) = \frac{1}{\epsilon^2}(H^n(\phi^{n+1} - \phi^n), U^{n+1}).
$$

Combining (3.34)-(3.40) and using the following identity

$$
(3.41) \quad (|\nabla\phi^n|^2, \rho^{n+1} - \rho^n) + (\rho^{n+1}, |\nabla\phi^{n+1}|^2 - |\nabla\phi^n|^2) = (\rho^{n+1}, |\nabla\phi^{n+1}|^2) - (\rho^n, |\nabla\phi^n|^2),
$$

we can obtain

$$
\begin{aligned}
(3.42) \quad &\frac{\beta}{2}\Big(\|\nabla\rho^{n+1}\|^2 - \|\nabla\rho^n\|^2 + \|\nabla(\rho^{n+1} - \rho^n)\|^2\Big) + \frac{1}{2}\Big(\|\nabla\phi^{n+1}\|^2 - \|\nabla\phi^{n+1}\|^2 + \|\nabla(\phi^{n+1} - \phi^n)\|^2\Big) \\
&+ \frac{\alpha}{2}\Big(\|\Delta\phi^{n+1}\|^2 - \|\Delta\phi^n\|^2 + \|\Delta(\phi^{n+1} - \phi^n)\|^2\Big) + \frac{1}{4\epsilon^2}\Big(\|U^{n+1}\|^2 - \|U^n\|^2 + \|U^{n+1} - U^n\|^2\Big) \\
&+ \frac{1}{4\eta^2}\Big(\|V^{n+1}\|^2 - \|V^n\|^2 + \|V^{n+1} - V^n\|^2\Big) - \theta\Big((\rho^{n+1}, |\nabla\phi^{n+1}|^2) - (\rho^n, |\nabla\phi^n|^2)\Big) \\
= &-M_\phi \delta t\|\nabla\mu_\phi^{n+1}\|^2 - M_\rho \delta t\|\nabla\mu_\rho^{n+1}\|^2.
\end{aligned}
$$

Finally, we obtain the desired result (3.32) after dropping some positive terms. □

**Remark 3.3.** *We note that the idea of the IEQ approach is very simple but quite different from the traditional time marching schemes. More precisely, it does not require the convexity required by the convex splitting approach (cf. [9]) or the boundedness for the second order derivative required by the stabilization approach (cf. [6,44,56]). Through a simple substitution of new variables, the complicated nonlinear potentials are transformed into quadratic forms. We summarize the great advantages of this quadratic transformations as follows: (i) this quadratization method works well for various complex nonlinear terms as long as the corresponding nonlinear potentials are bounded from below; (ii) the complicated nonlinear potential is transferred to a quadratic polynomial form which is much easier to handle; (iii) the derivative of the quadratic polynomial is linear, which provides the fundamental support for linearization method; (iv) the quadratic formulation in terms of new variables can automatically maintain this property of positivity (or bounded from below) of the nonlinear potentials.*

*When the nonlinear potential is the fourth order polynomial, e.g., the double well potential, the IEQ method is exactly the same as the so-called Lagrange Multiplier method developed in [15]. We remark that the idea of Lagrange Multiplier method only works well for the fourth order polynomial potential ($\phi^4$). This is because the nonlinear term $\phi^3$ (the derivative of $\phi^4$) can be naturally decomposed into a multiplication of two*

factors: $\lambda(\phi)\phi$ that is the Lagrange multiplier term, and the $\lambda(\phi) = \phi^2$ is then defined as the new auxiliary variable $U$. However, this method might not succeed for other type potentials. For instance, we notice the Flory-Huggins potential is widely used in two-phase model or fluid-surfactant model, see also [4, 53, 58]. The induced nonlinear term is logarithmic type as $\ln(\frac{\phi}{1-\phi})$. If one forcefully rewrites this term as $\lambda(\phi)\phi$, then $\lambda(\phi) = \frac{\ln(\frac{\phi}{1-\phi})}{\phi}$ that is the definition of the new variable $U$. Obviously, such a form is unworkable for algorithms design. About the application of the IEQ approach to handle other type of nonlinear potentials, we refer to the authors' other work in [58, 63, 71].

**Remark 3.4.** The IEQ approach provides more efficiency than the nonlinear approach. Let us consider the double well potential case, i.e., $E(\phi) = \int_\Omega (\phi^2 - 1)^2 d\boldsymbol{x}$, then IEQ scheme will generate the linear scheme as $(\phi^n)^2 \phi^{n+1}$. The implicit or convex splitting approach will produces the scheme as $(\phi^{n+1})^3$. Therefore, if the Newton iterative method is applied for this term, at each iteration the nonlinear convex splitting approach would yield the same linear operator as IEQ approach. Hence the cost of solving the IEQ scheme is the same as the cost of performing one iteration of Newton method for the implicit/convex splitting approach, provided that the same linear solvers are applied (for instance multi-grid with Gauss-Seidel relaxation). It is clear that the IEQ scheme would be much more efficient than the nonlinear schemes.

**Remark 3.5.** The proposed scheme follows the new energy dissipation law (3.11) formally instead of the energy law for the originated system (2.8). In the time-continuous case, the two energy laws are the same. In the time-discrete case, the energy $\mathbb{E}^{n+1}$ (defined in (3.33)) can be rewritten as a first order approximation to the Lyapunov functionals in $E(\phi^{n+1}, \rho^{n+1})$ (defined in (2.8)), that can be observed from the following facts heuristically. From (3.17), we have

$$U^{n+1} - ((\phi^{n+1})^2 - 1) = U^n - ((\phi^n)^2 - 1) + R_{n+1}, \tag{3.43}$$

where $R_{n+1} = O((\phi^{n+1} - \phi^n)^2)$. Since $R_k = O(\delta t^2)$ for $0 \leq k \leq n+1$ and $U^0 = (\phi^0)^2 - 1$, by mathematical induction we can easily get

$$U^{n+1} = (\phi^{n+1})^2 - 1 + O(\delta t). \tag{3.44}$$

**3.2. Second order Schemes.** For the IEQ approach, it is quite standard to develop second order version based on either BDF2 or Crank-Nicolson, the details are referred to the author's recent work in [58, 63, 71]. However, for the surfactant model studied in this paper, the nonlinear coupling term between $\rho$ and $\phi$ presents a huge barrier to obtain the energy stabilities. Even though the second order scheme is unconditionally energy stable for all numerical tests that we have performed in next section, it is still an open question about its stability proof or how to construct any second order linear schemes with unconditional energy stability. In this paper, following the strategy that we develop the first order scheme by combining the IEQ approach and explicit-implicit technique, we still list the second order BDF2 version as follows, for consistency.

<u>*Scheme*</u> **2.** (BDF2) Assuming that $\phi^{n-1}$, $\rho^{n-1}$, $U^{n-1}$, $V^{n-1}$ and $\phi^n$, $\rho^n$, $U^n$, $V^n$ are known, we solve $\phi^{n+1}$, $\rho^{n+1}$, $U^{n+1}$, $V^{n+1}$ from the following second order scheme:

**Step 1:** We update $\rho^{n+1}$ and $V^{n+1}$ as follows,

$$\frac{3\rho^{n+1} - 4\rho^n + \rho^{n-1}}{2\delta t} = M_\rho \Delta \mu_\rho^{n+1}, \tag{3.45}$$

$$\mu_\rho^{n+1} = -\beta \Delta \rho^{n+1} + \frac{1}{\eta^2} G^\star V^{n+1} - \theta |\nabla \phi^\star|^2, \tag{3.46}$$

$$3V^{n+1} - 4V^n + V^{n-1} = 2G^\star(3\rho^{n+1} - 4\rho^n + \rho^{n-1}). \tag{3.47}$$

**Step 2:** We update $\phi^{n+1}$ and $U^{n+1}$ as follows,

$$\frac{3\phi^{n+1} - 4\phi^n + \phi^{n-1}}{2\delta t} = M_\phi \Delta \mu_\phi^{n+1}, \tag{3.48}$$

$$\mu_\phi^{n+1} = -\Delta \phi^{n+1} + \alpha \Delta^2 \phi^{n+1} + \frac{1}{\epsilon^2} H^\star U^{n+1} + 2\theta \nabla \cdot (\rho^{n+1} \nabla \phi^{n+1}), \tag{3.49}$$

$$3U^{n+1} - 4U^n + U^{n-1} = 2H^\star(3\phi^{n+1} - 4\phi^n + \phi^{n-1}), \tag{3.50}$$

*where*

$$\phi^\star = 2\phi^n - \phi^{n-1}, H^\star = H(\phi^\star), G^\star = G(2\rho^n - \rho^{n-1}),$$

**Remark 3.6.** *The algorithm design for the hydrodynamics coupled model (2.10)-(2.13) presents further numerical challenges, e.g., how to decouple the computations of the velocity from the phase variables. This can be overcome by combining the proposed scheme in this paper for the two Cahn-Hilliard equations and the projection method [14] for the Navier-Stokes equations. Furthermore, one can use a decoupling technique that was developed in [?, 35, 46, 47, 69] to decouple the computations or $\phi$ and $\rho$ from the velocity field. Since this paper is focused on the development of efficient linear schemes for solving the nonlinearly coupled Cahn-Hilliard equations with multiple variables, the details of numerical schemes for the hydrodynamics coupled model that are in the similar vein as [28, 43, 45, 46, 73], are left to the interested readers.*

## 4. Numerical Simulations

We now present numerical experiments in two dimensions to validate the theoretical results derived in the previous section and demonstrate the efficiency, energy stability and accuracy of the proposed numerical schemes. The phase separation [24, 54] of dissimilar fluids are often carried out in the periodic domain, which is also adopted here. The computed domain is $[0, 2\pi]^2$, if not explicitly specified, the default values of order parameters are given below:

(4.1) $\quad M_\rho = 2.5\text{e}{-4}, M_\phi = 2.5\text{e}{-4}, \alpha = 2.5\text{e}{-4}, \beta = 1, \epsilon = 0.05, \eta = 0.08, \theta = 0.3, \rho_s = 1.$

We use the Fourier-spectral method to discretize the space, and $129^2$ Fourier modes are used so that the errors from the spatial discretization is negligible compared with the time discretization errors.

**Remark 4.1.** *The order parameter $\epsilon$ is chosen as the $0.8\% (= \frac{\epsilon}{2\pi})$ of the domain size, that is a reasonable choice for the computations of phase field models, see also [?, 6, 28–30, 40, 42–44, 46, 47, 59, 60, 69, 70, 73, 74]. In addition, note in quite a number of literatures, the choice of mobility parameter therein (denoted by $M_{old}$) is often chosen as 1, but our choice of the mobility parameter is $M = \epsilon^2 = 2.5\text{e}{-4}$ ($M_\rho$ and $M_\phi$) that seems to be much smaller. This is because the free energy in those literatures (denoted by $E_{old}$) is $E_{old} = \int_\Omega (\frac{\epsilon^2}{2}|\nabla\phi|^2 + (\phi^2-1)^2)d\boldsymbol{x}$ that is the rescale of the free energy given in this paper, namely, their Cahn-Hilliard equation reads as $\frac{\phi_t}{\epsilon^2 M_{old}} + \Delta(\Delta\phi + \frac{\phi(\phi^2-1)}{\epsilon^2}) = 0$, but our model reads as $\frac{\phi_t}{M} + \Delta(\Delta\phi + \frac{\phi(\phi^2-1)}{\epsilon^2}) = 0$. Therefore, the mobility parameter $M$ in this paper corresponds to $\epsilon^2 M_{old}$ instead of $M_{old}$ itself.*

### 4.1. Accuracy test.

We first test convergence rates of the two numerical schemes, the first order scheme (3.12)-(3.17) (denoted by LS1) and the second order BDF2 scheme (3.45)-(3.50) (denoted by LS2). The following initial conditions

(4.2) $\quad \begin{cases} \phi_0(x,y) = 0.3\cos(3x) + 0.5\cos(y), \\ \rho_0(x,y) = 0.2\sin(2x) + 0.25\sin(y) \end{cases}$

are used.

We perform the refinement test of the time step size, and choose the approximate solution obtained by using the scheme LS2 with the time step size $\delta t = 7.8125\text{e}{-5}$ as the benchmark solution (approximately the exact solution) for computing errors. We present the summations of the $L^2$ errors of two phase variables between the numerical solution and the exact solution at $t = 0.1$ with different time step sizes in Table 1. We observe that the schemes, LS1 and LS2, are first order and second order accurate respectively. Moreover, the second order scheme LS2 gives better accuracy than the first order scheme LS1 does when using the same time step. For instance, if the accuracy is required to be around $2\text{e}{-5}$, then for LS2, we can take the time step $5\text{e}{-3}$, but for LS1, we have to take the time step $6.25\text{e}{-4}$, therefore, the second order scheme LS2 can use the time step that is around $8 (= \frac{1\text{e}{-2}}{5\text{e}{-4}})$ times larger than the first order scheme LS1.

| $\delta t$ | LS1 | Order | LS2 | Order |
|---|---|---|---|---|
| 1e−2 | 4.21e−4 | – | 8.15e−5 | – |
| 5e−3 | 2.16e−4 | 0.96 | 2.18e−5 | 1.90 |
| 2.5e−3 | 1.09e−4 | 0.99 | 5.63e−6 | 1.95 |
| 1.25e−3 | 5.52e−5 | 0.98 | 1.42e−6 | 1.99 |
| 6.25e−4 | 2.77e−5 | 0.99 | 3.55e−7 | 2.00 |
| 3.125e−4 | 1.38e−5 | 1.00 | 8.48e−8 | 2.07 |
| 1.5625e−4 | 6.95e−6 | 0.99 | 2.10e−8 | 2.01 |

TABLE 1. The summation of $L^2$ numerical errors for $\phi$ and $\rho$ at $t = 0.1$ that are computed by schemes LS1, LS2 using various temporal resolutions with the initial conditions of (4.2), for mesh refinement test in time. The order parameters are of (4.1).

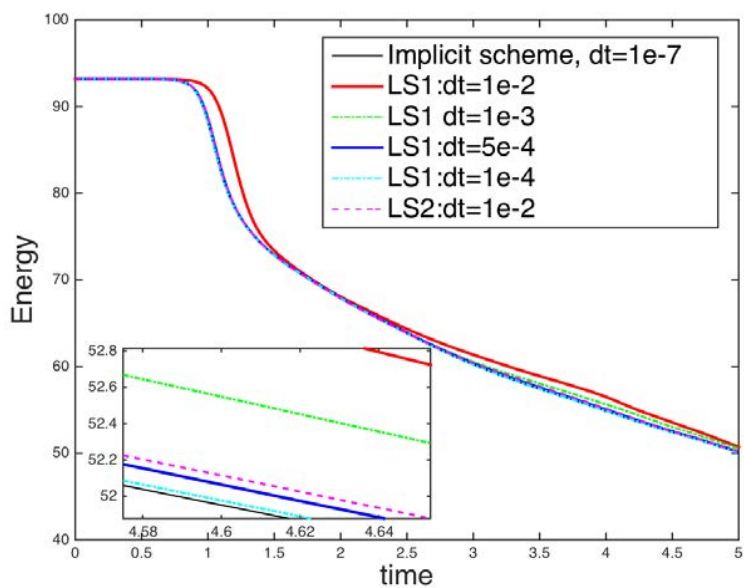

FIG. 1. Time evolution of the free energy functional for various time step sizes using the scheme LS1, LS2 and the second order fully implicit scheme. The small inset figure shows the small differences in the energy evolution for the considered time steps. The energy curves show the decays for all time step sizes, that confirms that our algorithm is unconditionally stable.

4.2. **Spinodal decomposition.** In this example, we study the phase separation dynamics that is called "spinodal decomposition". The process of the phase separation can be studied by considering a homogeneous binary mixture, which is quenched into the unstable part of its miscibility gap. In this case, the spinodal decomposition takes place, which manifests in the spontaneous growth of the concentration fluctuations that leads the system from the homogeneous to the two-phase state. Shortly after the phase separation starts, the domains of the binary components are formed and the interface between the two phases can be specified [2, 7, 75].

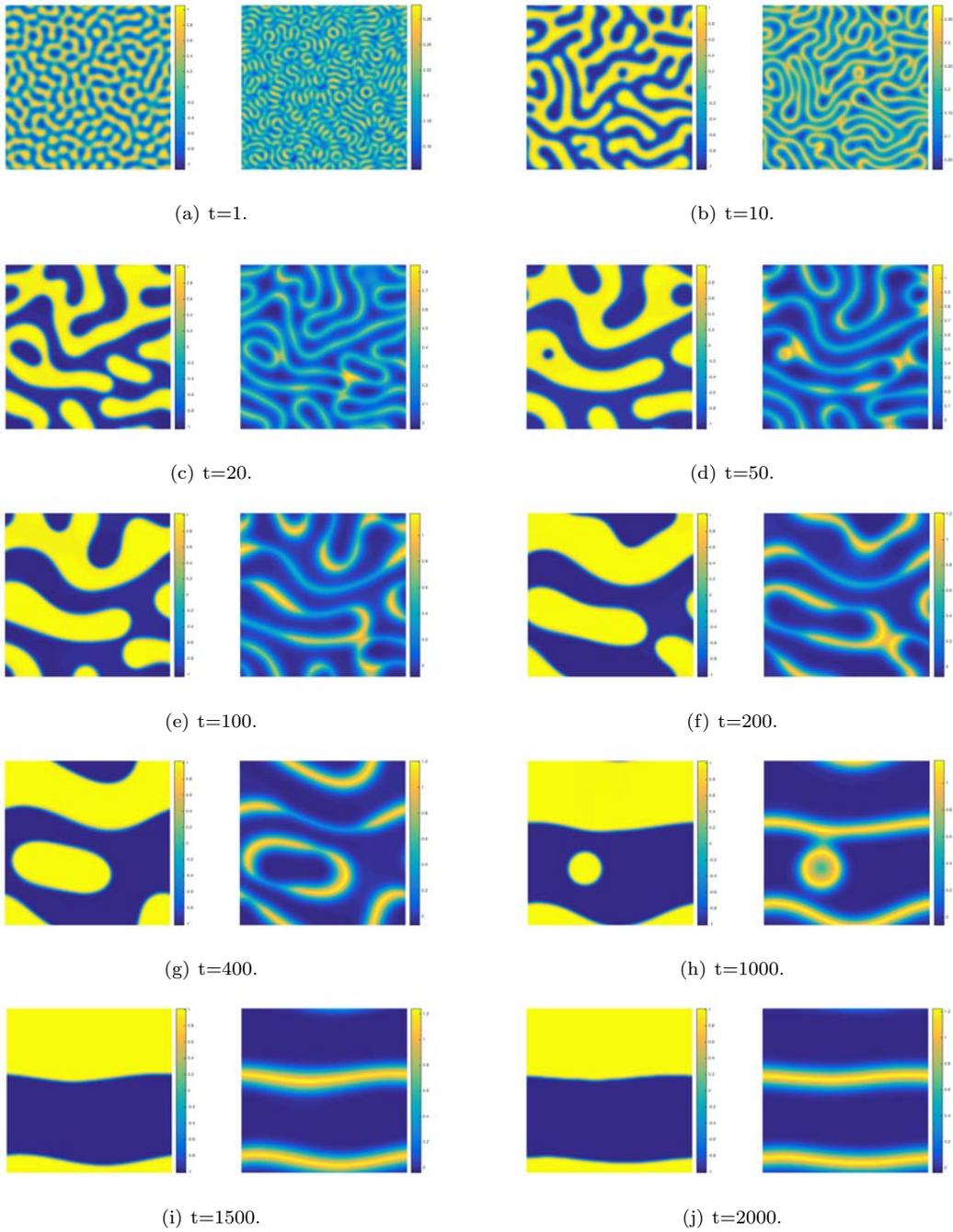

FIG. 2. Spinodal decomposition for $\overline{\phi}_0 = 0$ using parameters (4.1). Snapshots of phase variables $\phi$ and $\rho$ are taken at $t = 1, 10, 20, 50, 100, 200, 400, 1000, 1500, 2000$. For each subfigure, the left is the profile of $\phi$, and the right is the profile of $\rho$.

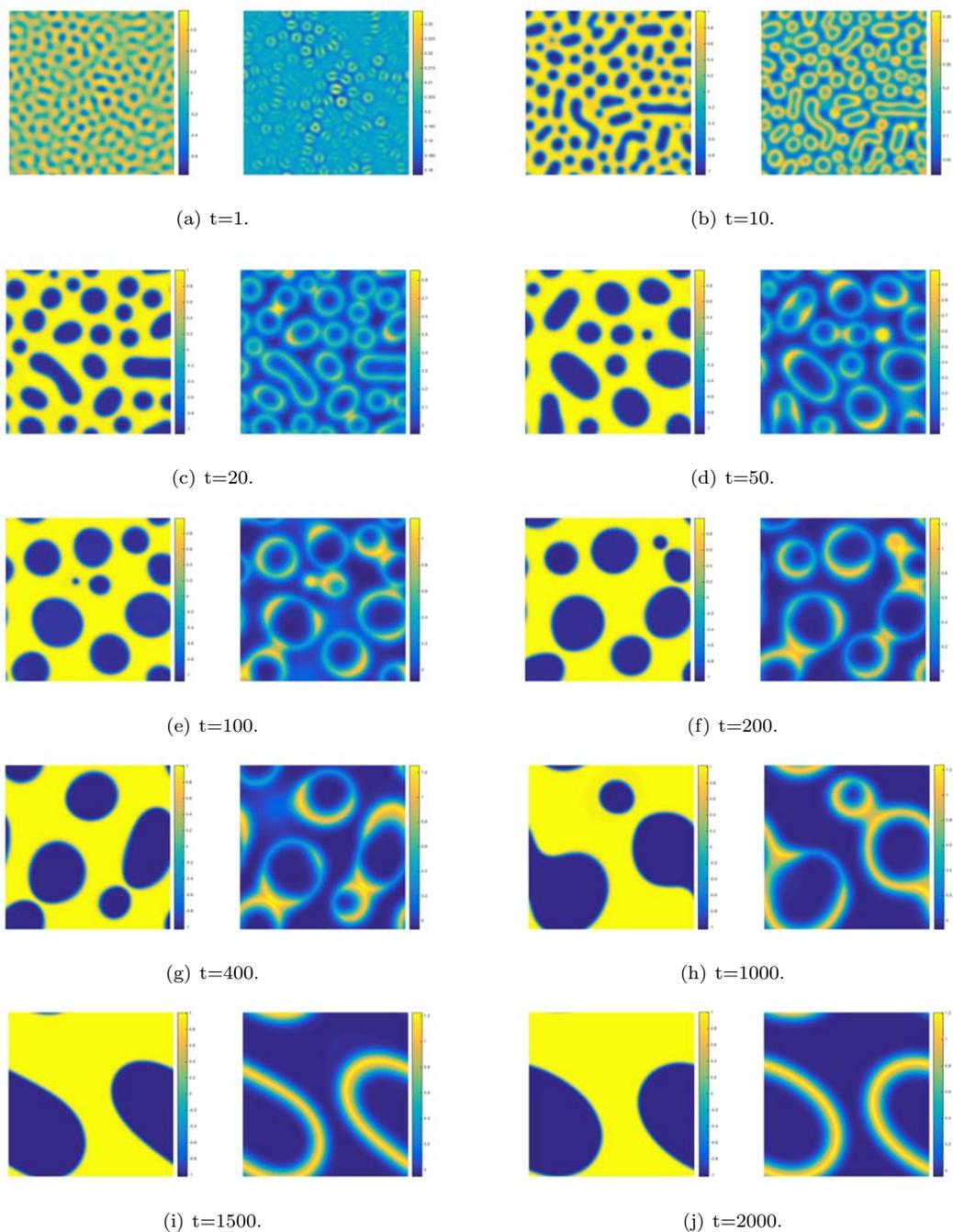

FIG. 3. Spinodal decomposition for $\overline{\phi}_0 = 0.2$ using parameters (4.1). Snapshots of phase variables $\phi$ and $\rho$ are taken at $t = 1, 10, 20, 50, 100, 200, 400, 1000, 1500, 2000$. For each subfigure, the left is the profile of $\phi$, and the right is the profile of $\rho$.

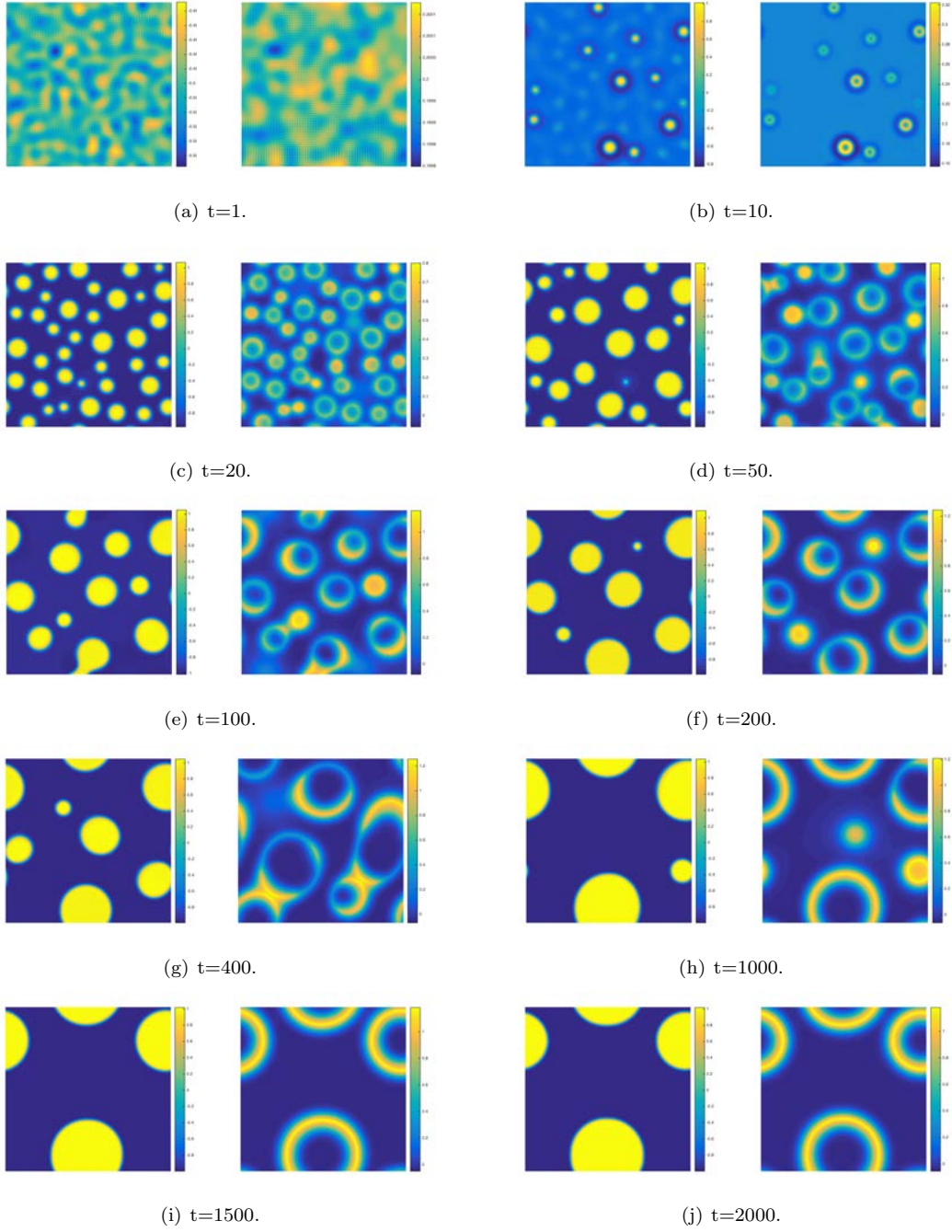

Fig. 4. Spinodal decomposition for $\overline{\phi}_0 = -0.5$ using parameters (4.1). Snapshots of phase variables $\phi$ and $\rho$ are taken at $t = 1, 10, 20, 50, 100, 200, 400, 1000, 1500, 2000$. For each subfigure, the left is the profile of $\phi$, and the right is the profile of $\rho$.

The initial conditions are taken as the randomly perturbed concentration fields as follows,

(4.3)
$$\phi_0(x,y) = \overline{\phi}_0 + 0.001 \mathrm{rand}(x,y),$$
$$\rho_0(x,y) = 0.2 + 0.001 \mathrm{rand}(x,y),$$

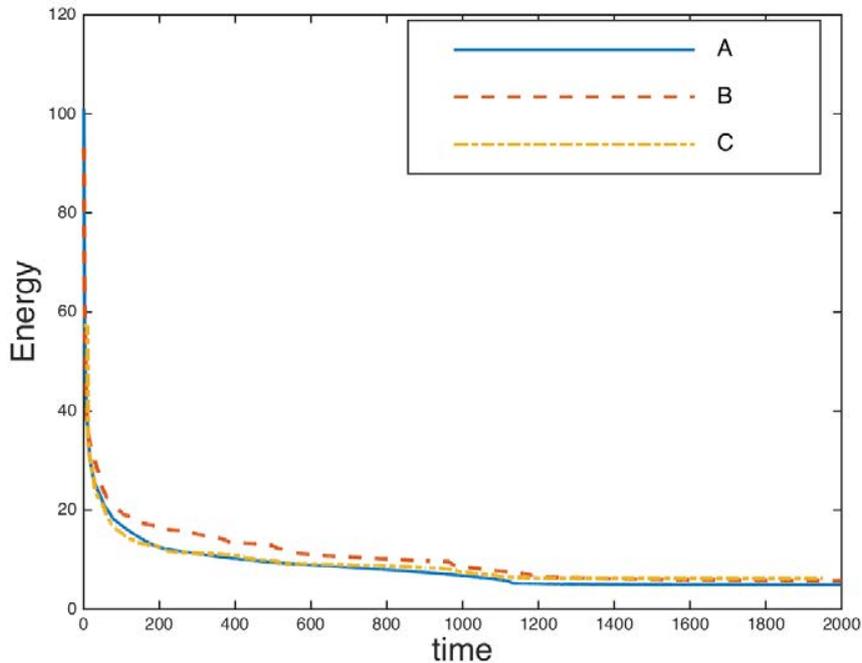

FIG. 5. Time evolution of the free energy functional for spinodal decomposition with A: $\overline{\phi}_0 = 0$, B: $\overline{\phi}_0 = 0.2$, C: $\overline{\phi}_0 = -0.5$ with order parameters (4.1).

where the rand$(x, y)$ is the random number in $[-1, 1]$ and has zero mean. We will vary $\overline{\phi}_0$ in next simulations.

We first choose $\overline{\phi}_0 = 0$ and compare the evolution of the free energy functional between the two proposed schemes and the second order fully implicit scheme with various time steps until $t = 5$ in Fig. 1. We take the computed solution of the fully implicit scheme by using a very tiny time step of $\delta t = 1e-7$ as the exact solution and compare with the results obtained from the schemes LS1 and LS2. For the first order scheme LS1, we choose four time steps with $\delta t = 1e-2, 1e-3, 5e-4$ and $1e-4$. We observe that all energy curves show the decays for all time step sizes, that confirm that our algorithms are unconditionally stable for any time step. Furthermore, when the time step size $\delta t$ is 0.01, the energy curve using the scheme LS1 is considerable (but not very far) away from the exact solution. For other smaller time steps, the obtained energy curves match well with the exact solution. This means the time step size has to be smaller than 0.01 at least in order to get reasonably good accuracy if using the scheme LS1 in practice. For the scheme LS2, we choose the time step $\delta t = 1e-2$, the obtained energy curve is already well matched with the exact solution. Moreover, we notice that its extent of agreement with the exact solution is even better than that of the scheme LS1 when using $\delta t = 5e-4$. In other words, under the same magnitude of accuracy requests, the second order scheme LS2 can use the time step that is around 20 ($= \frac{1e-2}{5e-4}$) times larger than the first order scheme LS1 at least. Therefore, we use the scheme LS2 and choose the time step $\delta t = 1e-3$ to perform the following spinodal decomposition simulations.

In Fig. 2, Fig. 3, and Fig. 4, we present the time evolution of the phase separation dynamics for $\overline{\phi}_0 = 0$, $0.2, -0.5$, respectively. Snapshots are taken at $t = 1, 10, 20, 50, 100, 200, 400, 1000, 1500$ and $2000$ for all cases. Initially, the two fluids are well mixed, and they sooner start to decompose due to surface tensions. A relatively high value of the concentration variable $\rho$ is always driven to be located at the fluid interface. For $\overline{\phi}_0 = 0$ in Fig. 2, the volume fractions of two fluids are the same, thus we observe the fluid interfaces are partially entangled and isolated. For other values of $\overline{\phi}_0 = 0.2, -0.5$, we observe that the fluids are accumulated to small drops finally. Overall, the numerical solutions in these simulations present similar features to those obtained in [21, 24–26, 54]. We plot the evolution of energy curves in Fig. 5 for all three cases, where the energy monotonically decays with respect to the time.

## 5. Concluding remarks

Although the phase field approach had been widely used to model the binary fluid-surfactant system for over twenty years, the development of energy stable schemes still remained very scarce due to its complex nature of nonlinearities. In this paper, we develop two linear, decoupled, computationally efficient schemes to solve the model, where the first order version is unconditionally energy stable. While we have considered only time discretizations here, the results can carry over to any consistent finite-dimensional Galerkin approximations (finite elements or spectral) since the proofs are all based on variational formulations with all test functions in the same space as the trial function. We also presented ample numerical results to validate the accuracy of the proposed schemes.